\documentclass[12pt]{article}
\usepackage{latexsym,amsfonts,amssymb}
\usepackage[dvips]{color}
\usepackage{colordvi,multicol}

\def \cal{\mathcal}

\newtheorem{thm}{Theorem}[section]
\newtheorem{cor}[thm]{Corollary}
\newtheorem{lem}[thm]{Lemma}

\date{}

\begin{document}
\title{\bf  The complement value problem for non-local
operators}
\author{}

\maketitle

\centerline{Wei Sun} \centerline{\small\it Department of
Mathematics and Statistics, Concordia University,}
\centerline{\small\it Montreal, H3G 1M8, Canada}
\centerline{\small\it  wei.sun@concordia.ca}

\vskip 1cm

%\vskip 1.4cm

\vskip 0.5cm \noindent Let $D$ be a bounded Lipschitz domain of
$\mathbb{R}^d$. We consider the complement value problem
$$
\left\{\begin{array}{l}(\Delta+a^{\alpha}\Delta^{\alpha/2}+b\cdot\nabla+c)u+f=0\ \ {\rm in}\ D,\\
u=g\ \ {\rm on}\ D^c.
\end{array}\right.$$
Under mild conditions, we show that there exists a unique bounded
continuous weak solution. Moreover, we give an explicit
probabilistic representation of the solution. The theory of
semi-Dirichlet forms and heat kernel estimates play an important
role in our approach.

\smallskip


\smallskip


\section{Introduction and Main Result}

Let $d\ge 1$ and $D$ be a bounded Lipschitz domain of
$\mathbb{R}^d$. Suppose $0<\alpha<2$ and $p>d/2$. Let $a>0$,
$b=(b_1,\dots,b_d)^*$ satisfying $|b|\in L^{2p}(D;dx)$ if $d\ge 2$
and $|b|\in L^{\infty}(D;dx)$ if $d=1$, $c\in L^{p\vee 2}(D;dx)$,
$f\in L^{2(p\vee 1)}(D;dx)$ and $g\in B_b(D^c)$. We consider the
complement value problem:
\begin{equation}\label{in1}
\left\{\begin{array}{l}(\Delta+a^{\alpha}\Delta^{\alpha/2}+b\cdot\nabla+c)u+f=0\ \ {\rm in}\ D,\\
u=g\ \ {\rm on}\ D^c.
\end{array}\right.
\end{equation}
The fractional Laplacian operator $\Delta^{\alpha/2}$ can be written in the form
$$
\Delta^{\alpha/2}\phi(x)=\lim_{\varepsilon\rightarrow 0}{\cal
A}(d,-\alpha)\int_{\{|x-y|\ge
\varepsilon\}}\frac{\phi(y)-\phi(x)}{|x-y|^{d+\alpha}}dy,\ \
\phi\in C^{\infty}_c(\mathbb{R}^d),
$$
where ${\cal A}(d,-\alpha):=\alpha2^{\alpha-1}\pi^{-d/2}\Gamma((d+\alpha)/2)\Gamma(1-\alpha/2)^{-1}$ and $C^{\infty}_c(\mathbb{R}^d)$ denotes the space of infinitely differentiable functions on $\mathbb{R}^d$ with compact
support.

The problem (\ref{in1}) is analogue of the Dirichlet problem for
second order elliptic integro-differential equations. For these
non-local equations, as opposed to the classical local case, the
function $g$ should be prescribed not only on the boundary
$\partial D$ but also in the whole complement $D^c$. The
complement value problem for non-local operators has many
applications, for example, in peridynamics \cite{AM,G,MD},
particle systems with long range interactions \cite{GL}, fluid
dynamics \cite{DG} and image processing \cite{GO}. The problem has
been widely studied by using different approaches from both
probability and analysis. These include, in particular, the
semi-group approach by Bony, Courr\`ege and Priouret \cite{BC},
the classical PDE approach by Garroni and Menaldi \cite{GM}, the
viscosity solution approach by Barles, Chasseigne and Imbert
\cite{B} and Arapostathisa, Biswasb and Caffarelli \cite{ABC}, and
the Hilbert space approach by Hoh and Jocob \cite{HJ} and
Felsinger, Kassmann and Voigt \cite{FKV}. Many results have also
been obtained for the interior and boundary regularity of
solutions, see for example, \cite{B,
B2,BL,Bo,BKK,CS1,CS2,CS,G2,G1,K,R1,R2} and the references therein.

 Different from \cite{ABC,B}, $b$,
$c$, $f$ and $g$ in (\ref{in1}) are not assumed to be continuous.
Also, the second order elliptic integro-differential operator in
(\ref{in1}) is not assumed to have the maximum principle. To
overcome these complications, in this paper, we will use the
theory of semi-Dirichlet forms to study both the existence and
uniqueness of solutions to the problem (\ref{in1}). Our work is
partially motivated by Guan and Ma \cite{MG}, which uses the
Dirichlet form approach to study the boundary value problem for
regional fractional Laplacians. The heat kernel estimates recently
obtained by Chen and Hu \cite{CH} play an important role in our
work.

Denote $L:=\Delta+a^{\alpha}\Delta^{\alpha/2}+b\cdot\nabla$. By
setting $b = 0$ off $D$, we may assume that the operator $L$ is
defined on $\mathbb{R}^d$. By \cite[Theorem 1.4]{CH}, the
martingale problem for $(L, C^{\infty}_c(\mathbb{R}^d))$ is
well-posed for every initial value $x\in \mathbb{R}^d$. We use
$((X_t)_{t\ge 0},({\cal F}_t)_{t\ge 0},(P_x)_{x\in\mathbb{R}^d})$ to denote the strong
Markov process associated with $L$. Let $\rho>0$.  Define
$$q_{\rho}(t,z)=t^{-d/2}\exp\left(-\frac{\rho|z|^2}{t}\right)+t^{-d/2}\wedge\frac{t}{|z|^{d+\alpha}},\
\ t>0,z\in\mathbb{R}^d.
$$ By
\cite[Theorems 1.2-1.4]{CH}, $X$ has a jointly continuous
transition density function $p(t,x,y)$ on $(0,\infty)\times
\mathbb{R}^d\times \mathbb{R}^d$, and for every $T>0$ there exist
positive constants $C_i,i=1,2,3,4$ such that
\begin{equation}\label{2}
C_1q_{C_2}(t,x-y)\le p(t,x,y)\le C_3 q_{C_4}(t,x-y),\ \ (t,x,y)\in
(0,T]\times\mathbb{R}^d\times\mathbb{R}^d.
\end{equation}

Define
$$
e(t):=e^{\int_0^tc(X_s)ds},\ \ t\ge0,
$$
and $\tau:=\inf\{t>0: X_t\in D^c\}$. Denote $W^{1,2}(D)=\{ u\in L^2(D;dx):|\nabla u|\in L^2(D;dx)\}$, $W^{1,2}_0(D)=\{ u\in W^{1,2}(D):\exists\{u_n\}_{n\in\mathbb{N}}\subset C_c^{\infty}(D)\ {\rm such\ that}\ u_n\rightarrow u\ {\rm in}\ W^{1,2}(D)\}$, and
$$
W^{1,2}_{loc}(D):=\{u: u\phi\in W^{1,2}_0(D)\ {\rm for\ any}\
\phi\in C^{\infty}_c(D)\}.
$$

The main result of this paper
is the following theorem.

\begin{thm}\label{thm1}
There exists $M>0$ such that if $\|c^+\|_{L^{p\vee 1}}\le M$, then
for any $f\in L^{2(p\vee 1)}(D;dx)$ and $g\in B_b(D^c)$, there
exists a unique $u\in B_b(\mathbb{R}^d)$ satisfying $u|_D\in W^{1,2}_{loc}(D)\cap
C(D)$ and
$$
\left\{\begin{array}{l}(\Delta+a^{\alpha}\Delta^{\alpha/2}+b\cdot\nabla+c)u+f=0\ \ {\rm in}\ D,\\
u=g\ \ {\rm on}\ D^c.
\end{array}\right.$$
Moreover, $u$ has the expression
\begin{equation}\label{31}
u(x)=E_x\left[e(\tau)g(X_{\tau})+\int_0^{\tau}e(s)f(X_s)ds\right],\
\ x\in \mathbb{R}^d.
\end{equation}
In addition, if $g$ is continuous at $z\in \partial D$ then
$$
\lim_{x\rightarrow z}u(x)=u(z).
$$
\end{thm}

Hereafter
 $(\Delta+a^{\alpha}\Delta^{\alpha/2}+b\cdot\nabla+c)u+f=0$ is
understood in the distribution sense: for any $\phi\in
C^{\infty}_c(D)$,
\begin{eqnarray}\label{1}
& &\int_{D}\langle \nabla u,\nabla\phi\rangle
dx+\frac{a^{\alpha}{\cal A}(d,-\alpha)}{2}\int_{\mathbb{R}^d}\int_{\mathbb{R}^d}\frac{(u(x)-u(y))(\phi(x)-\phi(y))}{|x-y|^{d+\alpha}}dxdy\nonumber\\
& &\ \ \ \ \ \ \ \ \ \ -\int_{D}\langle b,\nabla u\rangle \phi
dx-\int_{D}cu\phi dx-\int_{D}f\phi dx=0.
\end{eqnarray}
Note that the double integral appearing in (\ref{1}) is
well-defined for any $u\in B_b(\mathbb{R}^d)$ with $u|_D\in
W^{1,2}_{loc}(D)$ and $\phi\in C^{\infty}_c(D)$.

As a direct consequence of Theorem \ref{thm1}, we have the
following corollary.

\begin{cor}
If $c\le 0$, then for any $f\in L^{2(p\vee 1)}(D;dx)$ and $g\in
B_b(D^c)$ satisfying $g$ is continuous on $\partial D$, there
exists a unique $u\in B_b(\mathbb{R}^d)$ such that $u$ is
continuous on $\overline{D}$, $u|_D\in W^{1,2}_{loc}(D)$, and
$$
\left\{\begin{array}{l}(\Delta+a^{\alpha}\Delta^{\alpha/2}+b\cdot\nabla+c)u+f=0\ \ {\rm in}\ D,\\
u=g\ \ {\rm on}\ D^c.
\end{array}\right.$$
Moreover, $u$ has the expression
$$
u(x)=E_x\left[e(\tau)g(X_{\tau})+\int_0^{\tau}e(s)f(X_s)ds\right],\
\ x\in \mathbb{R}^d.
$$
\end{cor}

The proof of Theorem \ref{thm1} will be given in Section 3. In the
next section, we first present several lemmas. In particular, we will
use an old result of Kanda \cite{Kanda} to prove a key lemma (see
Lemma \ref{lem11} below), which will be used in proving the
continuity of solutions in Theorem \ref{thm1}.

\section{Some Lemmas}\setcounter{equation}{0}

 Throughout this paper, we denote by $(\cdot,\cdot)$ the inner product of $L^2(\mathbb{R}^d;dx)$ and denote by $C$ a generic fixed
strictly positive constant, whose value can change from line to
line. Recall that a measurable function $\varphi$ on $\mathbb{R}^d$ is said to be in the Kato class if and only if
$$
\left\{\begin{array}{ll}\lim\limits_{r\downarrow0}\left[\sup\limits_{x\in\mathbb{R}^d}\,\int_{\{|y-x|\le r\}}\frac{|\varphi(y)|}{|x-y|^{d-2}}dy\right]=0,\ \ &{\rm if}\ d\ge 3,\\
\lim\limits_{r\downarrow0}\left[\sup\limits_{x\in\mathbb{R}^d}\,\int_{\{|y-x|\le r\}}(-\ln(|x-y|))|\varphi(y)|dy\right]=0,\ \ &{\rm if}\ d= 2,\\
\sup\limits_{x\in\mathbb{R}^d}\,\int_{\{|y-x|\le 1\}}|\varphi(y)|dy<\infty,\ \ &{\rm if}\ d= 1.
\end{array}\right.
$$
\begin{lem}\label{l09}
Define
$$
\left\{\begin{array}{l}{\cal
E}^0(\phi,\psi)=\int_{\mathbb{R}^d}\langle
\nabla\phi,\nabla\psi\rangle
dx+\frac{a^{\alpha}{\cal A}(d,-\alpha)}{2}\int_{\mathbb{R}^d}\int_{\mathbb{R}^d}\frac{(\phi(x)-\phi(y))(\psi(x)-\psi(y))}{|x-y|^{d+\alpha}}dxdy\\
\ \ \ \ \ \ \ \ \ \ \ \ \ \ \ -\int_{\mathbb{R}^d}\langle
b,\nabla \phi\rangle \psi dx,\ \ \phi,\psi\in D({\cal E}^0),\\
D({\cal E}^0)=W^{1,2}(\mathbb{R}^d).
\end{array}
\right.
$$
Then, $({\cal E}^0,D({\cal E}^0))$ is a regular lower-bounded
semi-Dirichlet form on $L^2(\mathbb{R}^d;dx)$. Moreover,
$((X_t)_{t\ge 0},(P_x)_{x\in\mathbb{R}^d})$ is the Hunt process
associated with $({\cal E}^0,D({\cal E}^0))$.
\end{lem}

\noindent {\bf Proof.}\ \ By the assumption on $b$ and H\"older's inequality, we find that $|b|^2$ belongs to the Kato class. Then, we obtain by \cite[Chapter 7, Lemma 7.5]{Schechter} that there exists $\beta_0>0$ such that
\begin{equation}\label{beta0}
\int_{\mathbb{R}^d}|b|^2\phi^2dx\le\frac{1}{2}\int_{\mathbb{R}^d}|\nabla\phi|^2dx+\beta_0\int_{\mathbb{R}^d}|\phi|^2dx,\
\ \forall \phi\in W^{1,2}(\mathbb{R}^d).
\end{equation}
Define $$ {\cal E}^0_{\beta}(\phi,\psi)={\cal
E}^0(\phi,\psi)+\beta(\phi,\psi),\ \ \phi,\psi\in D({\cal E}^0).
$$
Then, $({\cal E}^0_{\beta},D({\cal E}^0))$ is a coercive closed
form on $L^2(\mathbb{R}^d;dx)$ for any $\beta>\beta_0$.

Denote by $C_0(\mathbb{R}^d)$ the space of continuous functions on $\mathbb{R}^d$ which vanish at infinity. If $\phi\in
C^{\infty}_c(\mathbb{R}^d)$, then $\Delta^{\alpha/2}\phi\in
C_0(\mathbb{R}^d)$ (cf. \cite[Theorem 31.5]{Sato}). Moreover, we
have $\Delta^{\alpha/2}\phi\in L^2(\mathbb{R}^d;dx)$. In fact,
suppose ${\rm supp}[\phi]\subset B(0,N)$ for some
$N\in\mathbb{N}$, then we get \begin{eqnarray*}
\int_{\{|x|>2N\}}|\Delta^{\alpha/2}\phi|^2dx&=&\int_{\{|x|>2N\}}\left({\cal A}(d,-\alpha)\int_{\mathbb{R}^d}\frac{\phi(x+y)}{|y|^{d+\alpha}}dy\right)^2dx\\
&=&\int_{\{|x|>2N\}}\left({\cal A}(d,-\alpha)\int_{\{|y|\ge1\}}\frac{\phi(x+y)}{|y|^{d+\alpha}}dy\right)^2dx\\
&\le&C\int_{\{|x|>2N\}}\int_{\{|y|\ge1\}}\frac{\phi^2(x+y)}{|y|^{d+\alpha}}dydx\\
&\le&C\int_{\mathbb{R}^d}\phi^2dx\int_{\{|y|\ge1\}}\frac{1}{|y|^{d+\alpha}}dy\\
&<&\infty.
\end{eqnarray*}
We have
$$
{\cal E}^0(\phi,\psi)=(-L\phi,\psi),\ \ \forall\phi,\psi\in
C^{\infty}_c(\mathbb{R}^d).
$$
By \cite[Theorem 3.1]{U}, $({\cal E}^0,D({\cal E}^0))$ is a
regular lower-bounded semi-Dirichlet form on
$L^2(\mathbb{R}^d;dx)$.

We now show that $((X_t)_{t\ge 0},(P_x)_{x\in\mathbb{R}^d})$ is
the Hunt process associated with $({\cal E}^0,D({\cal E}^0))$. We will follow the method of \cite[Section 4]{FU}, which relates the Hunt process associated with a semi-Dirichlet form to a martingale problem. Since $b$ in (\ref{in1}) is not assumed to be continuous, we cannot directly apply \cite[Theorem 4.3]{FU}. We give the detailed argument below.

Let
$((X^{\cal E}_t)_{t\ge 0},(P^{\cal E}_x)_{x\in\mathbb{R}^d})$ be a Hunt process
associated with $({\cal E}^0,D({\cal E}^0))$. Suppose that
$\phi\in C^{\infty}_c(\mathbb{R}^d)$. Define
$$
M^{\phi}_t:=\phi(X^{\cal E}_t)-\phi(X^{\cal
E}_0)-\int_0^tL\phi(X^{\cal E}_s)ds.
$$
Let $\psi=(1-L)\phi$. Then, $\psi\in L^2(\mathbb{R}^d;dx)$. Since
$\phi=G_1\psi$ $dx$-a.e., we get $\phi=R^{\cal E}_1\psi$ q.e.,
where $G_1$ and $R^{\cal E}_1$ are the 1-resolvents of ${\cal
E}^0$ and $X^{\cal E}$, respectively.  Hence
$$
M^{\phi}_t=R^{\cal E}_1\psi(X^{\cal E}_t)-R^{\cal E}_1\psi(X^{\cal
E}_0)-\int_0^t(R^{\cal E}_1\psi-\psi)(X^{\cal E}_s)ds, \ \ P^{\cal
E}_x-{\rm a.s.},\ {\rm q.e.}\ x\in \mathbb{R}^d,
$$
which implies that $\{M^{\phi}_t\}$ is a martingale under $P^{\cal
E}_x$ for q.e. $x\in \mathbb{R}^d$.

Let $\Phi$ be a countable subset of $C^{\infty}_c(\mathbb{R}^d)$
such that for any $\phi\in C^{\infty}_c(\mathbb{R}^d)$ there exist
$\{\phi_n\}\subset\Phi$ satisfying $\|\phi_n-\phi\|_{\infty},
\|\partial_i\phi_n-\partial_i\phi\|_{\infty},\|\partial_i\partial_j\phi_n-\partial_i\partial_j\phi\|_{\infty}\rightarrow0$
as $n\rightarrow\infty$ for any $i,j\in\{1,2,\dots,d\}$. Then, there is
an ${\cal E}^0$-exceptional set of $\mathbb{R}^d$, denoted by $F$,
such that $\{M^{\phi}_t\}$ is a martingale under $P^{\cal E}_x$
for any $x\in F^c$. Note that
$$
E^{\cal E}_x\left[\int_0^t|b\cdot\nabla\phi|(X^{\cal
E}_s)ds\right]\le e^t\||\nabla\phi|\|_{\infty}R^{\cal E}_1|b|(x).
$$
We obtain by taking limits that $\{M^{\phi}_t\}$ is a martingale
under $P^{\cal E}_x$ for any $\phi\in C^{\infty}_c(\mathbb{R}^d)$
and q.e. $x\in \mathbb{R}^d$. Therefore, by the uniqueness of
solutions to the martingale problem for $(L,
C^{\infty}_c(\mathbb{R}^d))$ (see \cite[Theorem 1.4]{CH}), we find
that $((X_t)_{t\ge 0},(P_x)_{x\in\mathbb{R}^d})$ is the Hunt
process associated with $({\cal E}^0,D({\cal E}^0))$.\hfill\fbox

\begin{lem}\label{lemw1}
$$\lim_{t\rightarrow 0}\sup_{x\in \mathbb{R}^d}P_x \left(\sup_{0\le
s\le t}|X_s-x|> r\right)=0, \ \forall r>0.
$$
\end{lem}

\noindent {\bf Proof.}\ \ Let $t,r>0$. Define
$$
\iota_t:=\sup_{x\in \mathbb{R}^d,0\le s\le t}P_x(|X_s-x|\ge
r)=\sup_{x\in \mathbb{R}^d,0\le s\le t}\int_{B(x,r)^c}p(s,x,y)dy.
$$
By (\ref{2}), we get \begin{equation}\label{sx1}
\lim_{t\rightarrow 0}\iota_t=0. \end{equation} Define
$$
S=\inf\{t>0:|X_t-X_0|>2r\}.
$$
For $x\in \mathbb{R}^d$, we have
\begin{eqnarray}\label{sx2}
P_x\left(\sup_{0\le s\le t}|X_s-x|> 2r\right)&=&P_x(S\le t)\nonumber\\
&\le&P_x(|X_t-x|\ge r)+P_x(S\le t, X_t\in B(x,r))\nonumber\\
&\le&\iota_t+P_x(S\le t\ {\rm and}\ |X_t-X_S|>r)\nonumber\\
&\le&\iota_t+E_x[1_{\{S\le t\}}P_{X_S}(|X_{t-S}-X_0|>r)]\nonumber\\
&\le&2\iota_t.
\end{eqnarray}
The proof is complete by (\ref{sx1}) and (\ref{sx2}).\hfill\fbox

Let $U$ be an open set of $\mathbb{R}^d$. Define
$$
\tau_U:=\inf\{t>0: X_t\in U^c\}.
$$
Denote by $p^U(t,x,y)$ the transition density function of the part
process $((X^U_t)_{t\ge 0},(P_x)_{x\in U})$. Define
$G_{\gamma}^U(x,y):=\int_0^{\infty}e^{-\gamma t}p^U(t,x,y)dt$ for
$x,y\in U$ and $\gamma\ge 0$.

\begin{lem}\label{add1} Let $U$ be a bounded open set of $\mathbb{R}^d$.

\noindent (1) For any $x\in U$,
\begin{equation}\label{0O}
P_x(\tau_U<\infty)=1.
\end{equation}

\noindent (2) There exist positive constants $\theta_1$ and $\theta_2$ such
that
\begin{equation}\label{3} p^U(t,x,y)\le \theta_1q_{\theta_2}(t,x-y),\ \
(t,x,y)\in (0,\infty)\times U\times U.
\end{equation}

\noindent (3) For any $t>0$, $P_x(\tau_U=t)=0$ and the function
$x\mapsto P_x(\tau_U>t)$ is upper semi-continuous on
$\mathbb{R}^d$.

\noindent (4) For any $x,y\in U$, the function $t\mapsto p^U(t,x,y)$ is continuous on $(0,\infty)$.

\end{lem}

\noindent {\bf Proof.}\ \  By (\ref{2}), similar to \cite[Lemma
6.1]{KS}, we can show that
\begin{equation}\label{0O1}
\sup_{x\in U}P_x(\tau_U>1)<1,
\end{equation}
and there exist positive constants $\theta^*_1$ and $\theta^*_2$
such that
\begin{equation}\label{0O2}
p^U(t,x,y)\le \theta^*_1e^{-\theta^*_2(t-1)},\ \ (t,x,y)\in
(1,\infty)\times U\times U.
\end{equation}
By (\ref{0O1}) and the Markov property of $X$, we conclude that
(\ref{0O}) holds. By (\ref{2}) and (\ref{0O2}), we conclude that
(\ref{3}) holds.

The proof of (3) is the same as \cite[Theorem 1.4.7 and Proposition 2.2.1]{PS}. We now prove (4). For $x,y\in U$ and $t>0$, we have
\begin{equation}\label{guan}
p^{U}(t,x,y)=p(t,x,y)-E_x[p(t-\tau_{U},X_{\tau_{U}},y)1_{\{\tau_{U}\le
t\}}].
\end{equation}
Then, (4) follows from (\ref{guan}), the continuity of $p(t,x,y)$, (\ref{2}) and (3).\hfill\fbox

\begin{lem}\label{Kato} Let $U$ be a bounded open set of $\mathbb{R}^d$. Suppose that $\varphi$ is a measurable function on $\mathbb{R}^d$ which belongs to the Kato class. Then, we have
$$
\lim_{t\rightarrow0}\sup_{x\in
U}E_x\left[\int_0^{t}|\varphi(X^U_s)|ds\right]=0.
$$

\end{lem}

\noindent {\bf Proof.}\ \ We have
\begin{equation}\label{4}
t^{-d/2}\wedge\frac{t}{|x-y|^{d+\alpha}}\le t^{-d/2}\le
e^{\rho}t^{-d/2}\exp\left(-\frac{\rho|x-y|^2}{t}\right)\ \ {\rm
if}\ |x-y|^2<t,
\end{equation}
and
\begin{eqnarray}\label{5}
& &\int_0^{|x-y|^2}\left(t^{-d/2}\wedge\frac{t}{|x-y|^{d+\alpha}}\right)dt\nonumber\\
&\le&\int_0^{|x-y|^2}\frac{t}{|x-y|^{d+\alpha}}dt\nonumber\\
&=&\frac{1}{2|x-y|^{d+\alpha-4}}.
\end{eqnarray}
Then, we obtain by (\ref{3}), (\ref{0O2}), (\ref{4}) and (\ref{5})
that there exists $C>0$ such that for any $x,y\in U$,
\begin{equation}\label{6}
G_0^U(x,y)\le\left\{\begin{array}{ll}\frac{C}{|x-y|^{d-2}}, & d\ge 3,\\
C\ln\left(1+\frac{1}{|x-y|}\right), &d=2,\\
C, & d=1.
\end{array}
\right.
\end{equation}
The proof is complete by Lemma \ref{lemw1}, (\ref{6}) and
\cite[Theorem 1]{Zhao}.\hfill\fbox

\begin{lem} There
exists $C>0$ such that
\begin{equation}\label{12}
\sup_{x\in D}E_x\left[\int_0^{\tau}v(X_s)ds\right]\le
C\|v\|_{L^{p\vee 1}},\ \ \forall v\in L_+^{p\vee 1}(D).
\end{equation}
\end{lem}

\noindent {\bf Proof.}\ \ We only prove (\ref{12}) when $d\ge 3$.
The cases that $d=1, 2$ can be considered similarly. Let $v\in
L_+^{p\vee 1}(D)$ and $x\in D$. Denote by $\varsigma(D)$ the diameter of $D$. By (\ref{6}), we have
\begin{eqnarray*}
E_x\left[\int_0^{\tau}v(X_s)ds\right]&\le&\int_DG_0^D(x,y)v(y)dy\\
&\le&C\int_D\frac{v(y)}{|x-y|^{d-2}}dy\\
&\le&C\left(\int_Dv(y)^pdy\right)^{1/p}\left(\int_D|x-y|^{-q(d-2)}dy\right)^{1/q}\\
 &=&C'\|v\|_{L^p}\left(\int_0^{\varsigma(D)}r^{d-1-q(d-2)}dr\right)^{1/q}\\
 &=&C^{''}\|v\|_{L^p},
\end{eqnarray*}
where $C'$ and $C^{''}$ are positive constants.  \hfill\fbox

Suppose that $\overline{D}\subset B(0,N)$ for some $N\in
\mathbb{N}$. Define \begin{equation}\label{Omega} \Omega=B(0,N).
\end{equation}
\begin{lem}\label{lem2}
Let $\gamma\ge 0$. For any compact set $K$ of $\Omega$, there
exist $\delta>0$ and $\vartheta_1,\vartheta_2\in (0,\infty)$ such
that for any $x,y\in K$ satisfying $|x-y|<\delta$, we have
\begin{equation}\label{shanghai1}
\left\{\begin{array}{ll}\frac{\vartheta_1}{|x-y|^{d-2}}\le G_{\gamma}^{\Omega}(x,y)\le \frac{\vartheta_2}{|x-y|^{d-2}},& {\rm if}\ d\ge3,\\
\vartheta_1\ln\frac{1}{|x-y|}\le G_{\gamma}^{\Omega}(x,y)\le \vartheta_2\ln\frac{1}{|x-y|},& {\rm if}\ d=2.
\end{array}\right.
\end{equation}
\end{lem}

\noindent {\bf Proof.}\ \ We only prove (\ref{shanghai1}) when
$d\ge 3$. The case that $d=2$ can be considered similarly.
Similar to (\ref{6}), we can prove that there exists
$\vartheta_2>0$ such that for any $x,y\in \Omega$,
$$
G_{\gamma}^{\Omega}(x,y)\le \frac{\vartheta_2}{|x-y|^{d-2}}.$$

We obtain by (\ref{2}) and (\ref{guan}) that there exist $C_1,C_2,C_3,\epsilon>0$ such that if $0<t\le \epsilon$ and $x,y\in K$
satisfying $|x-y|<\epsilon$ then
$$
p^{\Omega}(t,x,y)\ge
C_1t^{-d/2}\exp\left(-\frac{C_2|x-y|^2}{t}\right)-C_3.
$$
Thus, for $x,y\in K$
satisfying $|x-y|<\epsilon$, we have
\begin{eqnarray*}
G^{\Omega}_{\gamma}(x,y)&\ge&
e^{-\gamma\epsilon}\int_0^{\epsilon}p^{\Omega}(t,x,y)dt\nonumber\\
&\ge&e^{-\gamma\epsilon}\int_0^{\epsilon}\left[C_1t^{-d/2}\exp\left(-\frac{C_2|x-y|^2}{t}\right)-C_3\right]dt\nonumber\\
&\ge&e^{-\gamma\epsilon}\left[\int_0^{\infty}C_1t^{-d/2}\exp\left(-\frac{C_2|x-y|^2}{t}\right)dt-\int_{\epsilon}^{\infty}C_1t^{-d/2}dt-C_3\epsilon\right]\nonumber\\
&=& e^{-\gamma\epsilon}\left[\frac{C_4}{|x-y|^{d-2}}-\frac{C_5}{\epsilon^{\frac{d-2}{2}}}-C_3\epsilon\right],
\end{eqnarray*}
where $C_4$ and $C_5$ are positive constants. Therefore, there exist $0<\delta<\epsilon$ and $\vartheta_1>0$ such that if $x,y\in K$ satisfying
$|x-y|<\delta$ then
$$
G^{\Omega}_{\gamma}(x,y)\ge \frac{\vartheta_1}{|x-y|^{d-2}}.
$$\hfill\fbox

\begin{lem}\label{lem11}
Any point on $\partial D$ is a regular point of $D$ and $D^c$ for
the process $((X_t)_{t\ge 0},(P_x)_{x\in\mathbb{R}^d})$.
\end{lem}

\noindent {\bf Proof.}\ \ We first consider the case that $d\ge
2$. Let $\beta>\beta_0$ (see (\ref{beta0})) and $\Omega$ be
defined as in (\ref{Omega}). Denote by $((X^{\Omega}_t)_{t\ge 0},
(P^{\beta}_x)_{x\in \Omega})$ the Markov process associated with
$({\cal E}_{\beta}^0,W^{1,2}_0(\Omega))$. To prove the lemma, it
is sufficient to show that any point on $\partial D$ is a regular
point of $D$ and $D^c$ for the process $((X^{\Omega}_t)_{t\ge 0},
(P^{\beta}_x)_{x\in \Omega})$.

Let $A$ be a Borel set of $\Omega$ satisfying
$\overline{A}\subset\Omega$. Denote by $e_A$ the 0-equilibrium
measure of $A$ w.r.t $((X^{\Omega}_t)_{t\ge 0},
(P^{\beta}_x)_{x\in \Omega})$. Then, there exists a finite measure
$\mu_A$ concentrating on $\overline{A}$ such that (cf. \cite[page
58 and Theorem 3.5.1]{Oshima}),
$$
P^{\beta}_x(\sigma_A<\tau_{\Omega})=e_A(x)=\int_{\overline{A}}G^{\Omega}_{\beta}(x,y)\mu_A(dy)\
\ {\rm for\ q.e.}\ x\in\Omega,
$$
where $\sigma_A$ is the first hitting time of $A$. Since both
$\varphi(x):=P^{\beta}_x(\sigma_A<\tau_{\Omega})$ and
$\psi(x):=\int_{\overline{A}}G^{\Omega}_{\beta}(x,y)\mu_A(dy)$ are
0-excessive functions of $((X^{\Omega}_t)_{t\ge 0},
(P^{\beta}_x)_{x\in \Omega})$, we have
\begin{equation}\label{nan}
P^{\beta}_x(\sigma_A<\tau_{\Omega})=\int_{\overline{A}}G^{\Omega}_{\beta}(x,y)\mu_A(dy),\
\ \forall x\in\Omega.
\end{equation}

Let $z\in \partial D$. By the assumption on $D$, we know that $z$
is a regular point of $D$ and $D^c$ for the Brownian motion in
$\mathbb{R}^d$. Therefore, $z$ is a regular point of $D$ and $D^c$
for $((X_t)_{t\ge 0},(P_x)_{x\in\mathbb{R}^d})$ by Lemma
$\ref{lem2}$, (\ref{nan}) and \cite[Theorem 4.2]{Kanda}.

We now consider the case that $d=1$. To prove the lemma, it is
sufficient to show that for any $x\in \mathbb{R}^1$, $x$ is a
regular point of both $(-\infty,x)$ and $(x,\infty)$. We assume
without loss of generality that $x=0$. We will use an idea from
\cite{HSW} to show below that $0$ is a regular point of
$(0,\infty)$. Using the same method, we can show that $0$ is also
a regular point of $(-\infty, 0)$.

 Let $B$ be a Brownian motion
on $\mathbb{R}^1$ and $Y$ be a rotationally symmetric
$\alpha$-stable process on $\mathbb{R}^1$ that is independent of
$B$. Then, $B+aY$ is the symmetric L\'evy process associated with
$\Delta+a^{\alpha}\Delta^{\alpha/2}$. Denote by $\mathbb{P}$ and
$\mathbb{Q}$ the probability measures on
$D([0,\infty),\mathbb{R}^1)$ that are solutions to the martingale
problems for
$(\Delta+a^{\alpha}\Delta^{\alpha/2},C^{\infty}_c(\mathbb{R}^1))$
and $(L, C^{\infty}_c(\mathbb{R}^1))$ with initial value $0$,
respectively. Since $|b|\in L^{\infty}(D;dx)$, $\mathbb{P}$ and
$\mathbb{Q}$ are mutually locally absolutely continuous (cf. e.g.
\cite[Theorem 2.4]{CFY}). Define
$$
\sigma(\omega)=\inf\{t>0:\omega(t)>0\}, \ \ \sigma'(\omega)=\inf\{t>0:\omega(t)<0\} \ \ {\rm for}\ \omega\in D([0,\infty),\mathbb{R}^1),
$$
and
$$
S=\{\omega\in D([0,\infty),\mathbb{R}^1): \sigma(\omega)=0\}, \ \ S'=\{\omega\in D([0,\infty),\mathbb{R}^1): \sigma'(\omega)=0\}.
$$
By the Blumenthal 0-1 law, we know that $\mathbb{P}(S)=0$ or 1. If
$\mathbb{P}(S)=0$, then we obtain by the symmetry of $B+aY$ that
$\mathbb{P}(S')=0$ also. We have a contradiction. Therefore,
$$
\mathbb{P}(S)=1,
$$
which implies that
\begin{equation}\label{blu2}
\mathbb{P}(S^c)=0.\end{equation}

Define
$$
T_n=\{\omega\in D([0,\infty),\mathbb{R}^1): 0<\sigma(\omega)\le n\}\ \ {\rm for}\ n\in\mathbb{N},
$$
$$
T=\{\omega\in D([0,\infty),\mathbb{R}^1): 0<\sigma(\omega)<\infty\},
$$
$$
R_n=\{\omega\in D([0,\infty),\mathbb{R}^1): \sigma(\omega)>n\}\ \ {\rm for}\ n\in\mathbb{N},
$$
and
$$
R=\{\omega\in D([0,\infty),\mathbb{R}^1): \sigma(\omega)=\infty\}.
$$
Then, (\ref{blu2}) implies that $\mathbb{P}(T_n)=\mathbb{P}(R_n)=0$ for any $n\in\mathbb{N}$.
Since $\mathbb{Q}$ is locally absolutely continuous w.r.t. $\mathbb{P}$, we have $\mathbb{Q}(T_n)=\mathbb{Q}(R_n)=0$ for any $n\in\mathbb{N}$.
Then, $\mathbb{Q}(T)=\uparrow\mathbb{Q}(T_n)=0$ and $\mathbb{Q}(R)=\downarrow\mathbb{Q}(R_n)=0$. Therefore, $\mathbb{Q}(S)=1-\mathbb{Q}(T)-\mathbb{Q}(R)=1$,
which implies that $0$ is a regular point of $(0,\infty)$.\hfill\fbox

\begin{lem}\label{lemw2}
Define $\xi(x)=E_x[g(X_{\tau})]$ for $x\in\mathbb{R}^d$. If $g$ is
continuous at $z\in \partial D$, then $\lim_{x\rightarrow
z}\xi(x)=\xi(z)$.
\end{lem}

\noindent {\bf Proof.}\ \ Suppose that $g$ is continuous at $z\in
\partial D$. Let $\delta>0$. We define
$$
A_{\delta}=\{y\in \mathbb{R}^d: |y-z|<\delta\},\ \ \ \ \
T=\inf\{t>0: X_t\in A_{\delta}^c\}.
$$
For $t>0$, we have
$$\lim_{\stackrel{x\rightarrow z}{x\in D}}
P_x(T\le \tau)\le \limsup_{\stackrel{x\rightarrow z}{x\in
D}}P_x(\tau>t)+\limsup_{\stackrel{x\rightarrow z}{x\in D}}P_x(T\le
t).
$$
Then, we obtain by Lemma \ref{lemw1}, Lemma \ref{add1} (3)   and
Lemma \ref{lem11} that
\begin{equation}\label{tau}
\lim_{\stackrel{x\rightarrow z}{x\in D}}P_x(T\le \tau)=0.
\end{equation}

By the strong Markov property of $X$, we get
$$
\xi(x)=E_x[g(X_{\tau})1_{\{\tau< T\}}]+E_x[\xi(X_T)1_{\{\tau\ge
T\}}].
$$
Therefore, the proof is complete by the continuity of $g$ at $z$,
the boundedness of $g$ and (\ref{tau}).\hfill\fbox

\begin{lem}\label{lemp1}
For any $t>0$ and $z\in \partial D$, we have
\begin{equation}\label{upper}
\lim_{\stackrel{x\rightarrow z}{x\in D}}\left(\sup_{y\in
D}p^D(t,x,y)\right)=0.
\end{equation}
\end{lem}

\noindent {\bf Proof.}\ \ By (\ref{2}), for $\varepsilon<t$,  we
have
$$
p^D(t,x,y)=\int_Dp^D(\varepsilon,x,w)p^D(t-\varepsilon,w,y)dw\le
C(t-\varepsilon)^{-d/2}P_x(\tau>\varepsilon).
$$
Therefore, we obtain (\ref{upper}) by Lemma \ref{add1} (3) and Lemma \ref{lem11}.\hfill\fbox

\begin{lem}\label{lkjh}
Let $U$ be a bounded open set of $\mathbb{R}^d$ and
$\varphi\in B_b(\mathbb{R}^d)$ with ${\rm supp}[\varphi]\subset
\overline{U}^c$. Then, for $dx$-a.e. $x\in U$, we have
\begin{equation}\label{qa}
E_x[\varphi(X_{\tau_{U}})1_{\{\tau_{U}\le
t\}}]=a^{\alpha}{\cal A}(d,-\alpha)\int_0^t\left(\int_{\overline{U}^c}\int_{U}\frac{p^U(s,x,z)\varphi(y)}{|z-y|^{d+\alpha}}dzdy\right)ds.
\end{equation}
\end{lem}

\noindent {\bf Proof.}\ \ Let $\varphi\in B_b(\mathbb{R}^d)$ with
${\rm supp}[\varphi]\subset \overline{U}^c$ and $\psi\in
B_b(\mathbb{R}^d)$ with ${\rm supp}[\psi]\subset U$. By the
quasi-left continuity of $((X_t)_{t\ge
0},(P_x)_{x\in\mathbb{R}^d})$, we have
$$
E_{\psi\cdot dx}[\varphi(X_{\tau_{U}})1_{\{\tau_{U}\le
t\}}]=E_{\psi\cdot dx}[1_{\{X_{\tau_U-}\in
U\}}\varphi(X_{\tau_{U}})1_{\{\tau_{U}\le t\}}].
$$
By Lemma \ref{l09}, we know that $({\cal E}^0,W^{1,2}_0(U))$ is a
regular lower-bounded semi-Dirichlet form on $L^2(U;dx)$ and $X^U$
is the Hunt process associated with $({\cal E}^0,W^{1,2}_0(U))$
(cf. \cite[Theorem 3.5.7]{Oshima}). Let $({T}^U_t)_{t\ge
0}$ be the $L^2$-semigroup associated with $({\cal E}^0,W^{1,2}_0(U))$. Denote by $(\hat{T}^U_t)_{t\ge
0}$ the dual semigroup of $({T}^U_t)_{t\ge
0}$ on $L^2(U;dx)$. Similar to \cite[Lemma 4.5.5]{Fuku}, we can show that for any
$\varrho\in B_b(\mathbb{R}^d)$ with ${\rm supp}[\varrho]\subset
U$,
$$
E_{\psi\cdot
dx}[\varrho(X_{\tau_U-})\varphi(X_{\tau_{U}})1_{\{\tau_{U}\le
t\}}]=a^{\alpha}{\cal A}(d,-\alpha)\int_0^t\left(\int_{\overline{U}^c}\int_{U}\frac{\hat{T}^U_s\psi(x)\varrho(x)\varphi(y)}{|x-y|^{d+\alpha}}dxdy\right)ds.
$$
Then,
\begin{eqnarray*}
&&E_{\psi\cdot dx}[\varphi(X_{\tau_{U}})1_{\{\tau_{U}\le t\}}]\\
&=&a^{\alpha}{\cal A}(d,-\alpha)\int_0^t\left(\int_{\overline{U}^c}\int_{U}\frac{\hat{T}^U_s\psi(x)\varphi(y)}{|x-y|^{d+\alpha}}dxdy\right)ds\\
&=&a^{\alpha}{\cal A}(d,-\alpha)\int_{\mathbb{R}^d}\psi(x)\int_0^t\left(\int_{\overline{U}^c}\int_{U}\frac{p^U(s,x,z)\varphi(y)}{|z-y|^{d+\alpha}}dzdy\right)dsdx.
\end{eqnarray*}
Since $\psi$ is arbitrary, (\ref{qa}) holds for $dx$-a.e. $x\in
U$.\hfill\fbox

\section{Proof of Theorem \ref{thm1}}\setcounter{equation}{0}

\subsection{Boundedness and continuity of
solutions}\label{dfgh}

Let $u$ be defined by (\ref{31}). In this subsection, we will show
that $u\in B_b(\mathbb{R}^d)$, $u$ is continuous in $D$, and if
$g$ is continuous at $z\in
\partial D$ then $\lim_{x\rightarrow z}u(x)=u(z)$.

\noindent (1)\ \ By Khasminskii's inequality and (\ref{12}), there
exists $C>0$ such that for any $v\in L^{p\vee 1}_+(D)$ satisfying
$\|v\|_{L^{p\vee 1}}\le C$, we have
\begin{equation}\label{14}
\sup_{x\in D}E_x\left[e^{\int_0^{\tau}v(X_s)ds}\right]<\infty.
\end{equation}
In particular, this implies that there exists $\delta>0$ such that
\begin{equation}\label{15}
\sup_{x\in D}E_x\left[e^{\delta\tau}\right]<\infty.
\end{equation}
By (\ref{12}), we get
\begin{eqnarray}\label{13}
& &E_x\left[\left|\int_0^{\tau}e^{\int_0^sv(X_t)dt}f(X_s)ds\right|\right]\nonumber\\
&\le&\left(E_x\left[\int_0^{\tau}e^{\int_0^s2v(X_t)dt}ds\right]\right)^{1/2}\left(E_x\left[\int_0^{\tau}f^2(X_s)ds\right]\right)^{1/2}\nonumber\\
&\le&C\left(E_x\left[e^{\int_0^{\tau}2v(X_s)ds}\cdot\tau\right]\right)^{1/2}\|f^2\|^{1/2}_{L^{p\vee
1}}\nonumber\\
&\le&C\left(E_x\left[e^{\int_0^{\tau}4v(X_s)ds}\right]\right)^{1/4}\left(E_x\left[\tau^2\right]\right)^{1/4}\|f^2\|^{1/2}_{L^{p\vee
1}}.
\end{eqnarray}
By (\ref{14})--(\ref{13}), we know that there exists $M>0$ such
that if $\|c^+\|_{L^{p\vee 1}}\le M$, then for any $f\in
L^{2(p\vee 1)}(D;dx)$ and $g\in B_b(D^c)$, $u\in
B_b(\mathbb{R}^d)$.

\vskip 0.3cm \noindent (2)\ \ For $x\in D$ and $t>0$, we have
\begin{eqnarray}\label{re}
u(x)&=&E_x\left[e(\tau)g(X_{\tau})1_{\{\tau\le
t\}}+\int_0^{t\wedge
\tau}e(s)f(X_s)ds\right]\nonumber\\
& &+E_x\left[e(\tau)g(X_{\tau})1_{\{\tau> t\}}+1_{\{\tau>
t\}}\int_{t\wedge
\tau}^{\tau}e(s)f(X_s)ds\right]\nonumber\\
&=&E_x\left[e(\tau)g(X_{\tau})1_{\{\tau\le t\}}+\int_0^{t\wedge
\tau}e(s)f(X_s)ds\right]\nonumber\\
& &+E_x\left[e(t)1_{\{\tau>
t\}}E_{X_t}\left[e(\tau)g(X_{\tau})+\int_{0}^{\tau}e(s)f(X_s)ds\right]\right]\nonumber\\
&=&E_x\left[e(t)u(X_t)1_{\{\tau>
t\}}+e(\tau)g(X_{\tau})1_{\{\tau\le t\}}+\int_0^{t\wedge
\tau}e(s)f(X_s)ds\right].\ \
\end{eqnarray}
Define
$$
u_t(x)=E_x\left[u(X_t)\right],
$$
and
\begin{eqnarray*}\varepsilon_t(x)&=&E_x\left[-u(X_t)1_{\{\tau\le
t\}}+(e(t)-1)u(X_t)1_{\{\tau> t\}}+e(\tau)g(X_{\tau})1_{\{\tau\le
t\}}+\int_0^{t\wedge \tau}e(s)f(X_s)ds\right]\\
&:=&\sum_{i=1}^4\varepsilon^{(i)}_t.
\end{eqnarray*}
Then, we have $u=u_t+ \varepsilon_t$. By (\ref{2}) and the joint
continuity of $p(t,x,y)$ on $(0,\infty)\times \mathbb{R}^d\times
\mathbb{R}^d$, we obtain that $u_t$ is continuous in $D$. By Lemma
\ref{lemw1}, we find that
\begin{equation}\label{d1}
\lim_{t\rightarrow 0}P_x(\tau\le t)=0\ \ {\rm uniformly\ on\ any\
compact\ subset\ of}\ D.
\end{equation}
Then, we obtain by the boundedness of $u$ and (\ref{d1}) that
$\varepsilon^{(1)}_t$ converges to 0 uniformly on any compact
subset of $D$.

Let $\varphi=|c|+|f|$. By Lemma \ref{Kato} and the assumptions on $c$ and $f$, we have
\begin{equation}\label{a1}
\lim_{t\rightarrow 0}\sup_{x\in
D}E_x\left[\int_0^{t}\varphi(X^D_s)ds\right]=0,
\end{equation}
and
\begin{equation}\label{a2}\sup_{x\in
D}E_x\left[e^{\int_0^{t}r\varphi(X^D_s)ds}\right]<\infty,\ \ \forall
t,r>0.
\end{equation}
Note that for $t<\tau$, we have
\begin{eqnarray}\label{fitz1}
e(t)&=&1-\left.\left(e^{\int_s^tc(X^D_w)dw}\right)\right|_{s=0}^t\nonumber\\
&=&1-\int_0^td\left(e^{\int_s^tc(X^D_w)dw}\right)\nonumber\\
&=&1+\int_0^te^{\int_s^tc(X^D_w)dw}c(X^D_s)ds.
\end{eqnarray}
By (\ref{a1}), (\ref{a2}) and (\ref{fitz1}), we get
\begin{eqnarray*}
&&\lim_{t\rightarrow 0}\sup_{x\in D}E_x[1_{\{\tau>
t\}}|e(t)-1|]\nonumber\\
&\le&\lim_{t\rightarrow 0}\sup_{x\in D}E_x\left[\int_0^{t}e^{\int_s^{t}|c|(X^D_w)dw}|c|(X^D_s)ds\right]\nonumber\\
&=&\lim_{t\rightarrow 0}\sup_{x\in
D}E_x\left[\int_0^{t}|c|(X^D_s)E_{X^D_s}\left[e^{\int_0^{t-s}|c|(X^D_w)dw}\right]ds\right]\nonumber\\
&=&0.
\end{eqnarray*}
Then, $\varepsilon^{(2)}_t$ converges to 0 uniformly on $D$.

By (\ref{d1}), (\ref{a2}) and the boundedness of $g$, we obtain
that $\varepsilon^{(3)}_t$ converges to 0 uniformly on any compact
subset of $D$.  Similar to (\ref{fitz1}), we can show that for $t<\tau$,
\begin{equation}\label{fitz2}
e^{\int_0^t\varphi(X^D_s)ds}=1+\int_0^te^{\int_s^t\varphi(X^D_w)dw}\varphi(X^D_s)ds.
\end{equation}
By (\ref{a1}), (\ref{a2}) and (\ref{fitz2}), we get
\begin{eqnarray*}
&&\lim_{t\rightarrow 0}\sup_{x\in D}\left|E_x\left[\int_0^{t\wedge \tau}e(s)f(X_s)ds\right]\right|\nonumber\\
&\le&\lim_{t\rightarrow 0}\sup_{x\in D}E_x\left[\int_0^{t}e^{\int_0^s\varphi(X^D_w)dw}\varphi(X^D_s)ds\right]\nonumber\\
&=&\lim_{t\rightarrow 0}\sup_{x\in D}E_x\left[\int_0^{t}e^{\int_s^{t}\varphi(X^D_w)dw}\varphi(X^D_s)ds\right]\nonumber\\
&=&\lim_{t\rightarrow 0}\sup_{x\in
D}E_x\left[\int_0^{t}\varphi(X^D_s)E_{X^D_s}\left[e^{\int_0^{t-s}\varphi(X^D_w)dw}\right]ds\right]\nonumber\\
&=&0.
\end{eqnarray*}
Then,
$\varepsilon^{(4)}_t$ converges to 0 uniformly on $D$. Therefore,
$u$ is continuous in $D$.

\vskip 0.3cm \noindent (3)\ \  Define
\begin{equation}\label{g2}
{\cal M}_t=u(X_t)1_{\{\tau> t\}}+g(X_{\tau})1_{\{\tau\le
t\}}+\int_0^{t\wedge \tau}(f+cu)(X_s)ds,
\end{equation}
and
\begin{equation}\label{h1}
{\cal N}_t=e(t)u(X_t)1_{\{\tau>t\}}+e(\tau)g(X_{\tau})1_{\{\tau\le
t\}}+\int_0^{t\wedge\tau}e(s)f(X_s)ds.
\end{equation}

Let $0\le s<t$. By (\ref{re}), we get
$$
u(X_s)=E_{X_s}\left[e(t-s)u(X_{t-s})1_{\{\tau>
t-s\}}+e(\tau)g(X_{\tau})1_{\{\tau\le t-s\}}+\int_0^{(t-s)\wedge
\tau}e(w)f(X_w)dw\right],
$$
which together with the strong Markov property of $X$ implies that
$$
E_x[{\cal N}_t-{\cal N}_s|{\cal F}_s]=0.
$$
Then, $({\cal N}_t)_{t\ge 0}$ is a
martingale under $P_x$ for any $x\in D$.

By (\ref{g2}) and (\ref{h1}), we get \begin{eqnarray}\label{dddd}
{\cal N}_t&=&e(t)u(X_t)1_{\{\tau>t\}}+e(t)g(X_{\tau})1_{\{\tau\le
t\}}-\int_0^te(s)c(X_s)g(X_{\tau})1_{\{\tau\le s\}}ds\nonumber\\
& &+\int_0^{t\wedge\tau}e(s)f(X_s)ds\nonumber\\
&=&e(t)u(X_t)1_{\{\tau>t\}}+e(t)g(X_{\tau})1_{\{\tau\le t\}}+e(t)\int_0^{t\wedge \tau}(f+cu)(X_s)ds\nonumber\\
& &-\int_0^te(s)c(X_s)\left(u(X_s)1_{\{\tau>
s\}}+g(X_{\tau})1_{\{\tau\le
s\}}+\int_0^sf(X_w)1_{\{\tau\ge w\}}dw\right)ds\nonumber\\
& &-\int_0^tc(X_w)u(X_w)1_{\{\tau\ge w\}}\left(\int_w^t
e(s)c(X_s)ds\right)dw\nonumber\\
&=&e(t){\cal M}_t-\int_0^t{\cal M}_sde(s).
\end{eqnarray}
By the integration by parts formula for semi-martingales, we have
$$
e(t){\cal M}_t-u(x)=\int_0^t{\cal M}_sde(s)+\int_0^te(s)d{\cal
M}_s.
$$
Hence we obtain by (\ref{dddd}) that $({\cal M}_t)_{t\ge 0}$ is a
martingale under $P_x$ for any $x\in D$. Therefore, we have
\begin{equation}\label{g4}
u(x)=E_x\left[u(X_t)1_{\{\tau> t\}}+g(X_{\tau})1_{\{\tau\le
t\}}+\int_0^{t\wedge \tau}(f+cu)(X_s)ds\right], \ \ x\in D.
\end{equation}

Define \begin{equation}\label{xi} \xi(x)=E_x[g(X_{\tau})],\ \ x\in
\mathbb{R}^d,\end{equation}
 and \begin{equation}\label{w}
 w(x)=u(x)-\xi(x), \ \ x\in
\mathbb{R}^d. \end{equation} By (\ref{g4}), we get
\begin{equation}\label{xw}
w(x)=E_x\left[w(X_t)1_{\{\tau>t\}}+\int_0^{t\wedge\tau}(f+cu)(X_s)ds\right],
\ \ x\in D.
\end{equation}
By the assumptions on $f$ and $c$, the boundedness of $u$ and
Lemma \ref{Kato}, we have
\begin{equation}\label{a11}
\lim_{t\rightarrow 0}\sup_{x\in
D}E_x\left[\int_0^{t\wedge\tau}|f+cu|(X_s)ds\right]=0.
\end{equation}
Therefore, we obtain by Lemma \ref{lemw2}, Lemma \ref{lemp1} and
(\ref{xi})--(\ref{a11}) that if $g$ is continuous at $z\in
\partial D$, then $\lim_{x\rightarrow z}u(x)=u(z)$.

\subsection{Existence of solutions}

Let $u$ be defined by (\ref{31}), and $\xi$ and $w$ be defined by
(\ref{xi})
 and (\ref{w}), respectively.

We will first show that $\xi\in
 W^{1,2}_{loc}(D)$ and ${\cal E}^0(\xi,\phi)=0$ for
any $\phi\in C^{\infty}_c(D)$. We assume without loss of
generality that $g\ge 0$ on $D^c$. Let $\{D_n\}_{n\in\mathbb{N}}$
be a sequence of relatively compact open subsets of $D$ such that
$\overline{D}_n\subset D_{n+1}$ and $D=\cup_{n=1}^{\infty}D_n$,
and $\{\chi_n\}_{n\in\mathbb{N}}$  be a sequence of functions in
$C^{\infty}_c(D)$ such that $0\le \chi_n\le 1$ and
$\chi_n|_{D_n}=1$. Suppose that $\beta>\beta_0$ (see (\ref{beta0})).
Let $e^{\beta}_{D_n}$ be the $\beta$-equilibrium of $D_n$ w.r.t.
$X^D$. By \cite[Lemma 2.1.1]{Oshima}, $e^{\beta}_{D_n}\in
W^{1,2}_0(D)$ and $e^{\beta}_{D_n}=1$ $dx$-a.e. on $D_n$. Note
that
\begin{equation}\label{xi11}
\xi(x)=E_x[\xi(X_\tau)], \ \ x\in \mathbb{R}^d. \end{equation} We
find that $\xi|_D$ is a $\beta$-excessive function w.r.t. $X^D$.
Then, we get $(\|\xi\|_{\infty}e^{\beta}_{D_n})\wedge \xi\in
W^{1,2}_0(D)$ (cf. \cite[Theorem 2.6]{MOR}). Since
$(\|\xi\|_{\infty}e^{\beta}_{D_n})\wedge \xi=\xi$ $dx$-a.e. on
$D_n$ and $n\in\mathbb{N}$ is arbitrary, we have $\xi\in
W^{1,2}_{loc}(D)$.

Suppose $\phi\in C^{\infty}_c(D_m)$ for some $m\in\mathbb{N}$. By
(\ref{xi11}), we know that $(\xi(X_{t\wedge\tau}))_{t\ge 0}$ is a
martingale under $P_x$ for $x\in D$. By the integration by parts
formula for semi-martingales, we get
$$
E_x[e^{-\beta(t\wedge\tau_{D_m})}\xi(X_{t\wedge\tau_{D_m}})]=\xi(x)-\beta E_x\left[\int_0^{t\wedge\tau_{D_m}}e^{-\beta
s}\xi(X_s)ds\right].
$$
Then, we have
\begin{equation}\label{ch1}
\lim_{t\rightarrow
0}\int_{\mathbb{R}^d}\phi(x)\frac{\xi(x)-E_x[e^{-\beta(t\wedge\tau_{D_m})}\xi(X_{t\wedge\tau_{D_m}})]}{t}dx=\beta\int_{D_m}\xi\phi
dx.
\end{equation}
For $n>m$, define
$$
\eta_n(x)=E_x[e^{-\beta{\tau_{D_m}}}(\xi\chi_n)(X_{\tau_{D_m}})],\
\ x\in\mathbb{R}^d.
$$
We have
$\eta_n(x)=E_x[e^{-\beta{(t\wedge\tau_{D_m}})}\eta_n(X_{t\wedge\tau_{D_m}})]$
for $t\ge 0$ and $x\in D_m$, and $\eta_n(x)=\xi\chi_n(x)$ for
q.e.-$x\in D^c_m$. By \cite[Theorem 3.5.1]{Oshima}, we get
\begin{eqnarray}\label{ch2}
&&{\cal E}^0_{\beta}(\xi\chi_n,\phi)\nonumber\\
&=&{\cal
E}^0_{\beta}(\xi\chi_n-\eta_n,\phi)\nonumber\\
&=&\lim_{t\rightarrow
0}\int_{D_m}\phi(x)\frac{\xi\chi_n-\eta_n-E_x[e^{-\beta(t\wedge\tau_{D_m})}((\xi\chi_n)(X_{t\wedge\tau_{D_m}})-\eta_n(X_{t\wedge\tau_{D_m}}))]}{t}dx\nonumber\\
&=&\lim_{t\rightarrow
0}\int_{D_m}\phi(x)\frac{\xi-E_x[e^{-\beta(t\wedge\tau_{D_m})}(\xi\chi_n)(X_{t\wedge\tau_{D_m}})]}{t}dx.
\end{eqnarray}

By (\ref{ch1}) and (\ref{ch2}), we get
\begin{eqnarray}\label{fitz3}
&&{\cal E}^0(\xi\chi_n,\phi)\nonumber\\
&=&\lim_{t\rightarrow
0}\frac{1}{t}\int_{D_m}\phi(x){E_x[e^{-\beta(t\wedge\tau_{D_m})}\xi(X_{t\wedge\tau_{D_m}})(1-\chi_n(X_{t\wedge\tau_{D_m}}))]}dx\nonumber\\
&=&\lim_{t\rightarrow
0}\frac{1}{t}\int_{D_m}\phi(x){E_x[1_{\{\tau_{D_m}\le t\}}e^{-\beta\tau_{D_m}}\xi(X_{t\wedge\tau_{D_m}})(1-\chi_n(X_{t\wedge\tau_{D_m}}))]}dx.
\end{eqnarray}
Let $\varepsilon>0$. There exists $\delta>0$ such that for any $0<t<\delta$, $1-e^{-\beta t}<\varepsilon$. Suppose that $\overline{D}\subset B(0,N)$ for some $N\in
\mathbb{N}$. Then, we obtain by Lemma \ref{lkjh} that for $0<t<\delta$,
\begin{eqnarray}\label{fitz4}
&&\left|\frac{1}{t}\int_{D_m}\phi(x){E_x[1_{\{\tau_{D_m}\le t\}}(1-e^{-\beta\tau_{D_m}})\xi(X_{t\wedge\tau_{D_m}})(1-\chi_n(X_{t\wedge\tau_{D_m}}))]}dx\right|\nonumber\\
&\le&\frac{\varepsilon}{t}\int_{D_m}|\phi(x)|{E_x[1_{\{\tau_{D_m}\le t\}}\xi(X_{t\wedge\tau_{D_m}})(1-\chi_n(X_{t\wedge\tau_{D_m}}))]}dx\nonumber\\
&=&\frac{\varepsilon a^{\alpha}{\cal A}(d,-\alpha)}{t}\int_{D_m}|\phi(x)|\left[\int_0^t\int_{\overline{D}^c_m}\left(\xi(y)(1-\chi_n(y))\int_{D_m}\frac{p^{D_m}(s,x,z)}{|z-y|^{d+\alpha}}dz\right)dyds\right]dx\nonumber\\
&\le&\frac{\varepsilon a^{\alpha}{\cal A}(d,-\alpha)}{t}\left(\int_{D_m}|\phi(x)|\left[\int_0^t\int_{(B(0,2N))^c}\left(\xi(y)
(1-\chi_n(y))\int_{D_m}\frac{p^{D_m}(s,x,z)}{|z-y|^{d+\alpha}}dz\right)dyds\right]dx\right.\nonumber\\
&& +\left.\int_{D_m}|\phi(x)|\left[\int_0^t\int_{B(0,2N)\cap D^c_n}\left(\xi(y)
(1-\chi_n(y))\int_{D_m}\frac{p^{D_m}(s,x,z)}{|z-y|^{d+\alpha}}dz\right)dyds\right]dx\right)\nonumber\\
&\le&\varepsilon a^{\alpha}{\cal A}(d,-\alpha)\|\phi\|_{\infty}\|\xi\|_{\infty}|D_m|\left(\int_{(B(0,2N))^c}\frac{1}{(|y|/2)^{d+\alpha}}dy+\vartheta^{-(d+\alpha)}|B(0,2N)\cap
D^c_n|\right),
\end{eqnarray}
where $\vartheta=\inf\{|x-y|:x\in D_m, y\in D_n^c\}$, and $|D_m|$ and
$|B(0,2N)\cap D^c_n|$ denote the Lebesgue measures of $D_m$ and $B(0,2N)\cap
D_n^c$, respectively. Since $\varepsilon>0$ is arbitrary, we obtain by (\ref{fitz3}) and (\ref{fitz4}) that
\begin{equation}\label{fitz5}
{\cal E}^0(\xi\chi_n,\phi)=\lim_{t\rightarrow
0}\frac{1}{t}\int_{D_m}\phi(x){E_x[1_{\{\tau_{D_m}\le t\}}\xi(X_{t\wedge\tau_{D_m}})(1-\chi_n(X_{t\wedge\tau_{D_m}}))]}dx.
\end{equation}

Define
\begin{equation}\label{weis1}
F_n(z)=\int_{D^c_n}\frac{\xi(y)(1-\chi_n(y))}{|z-y|^{d+\alpha}}dy,\ \ z\in D_m.
\end{equation}
Then, $F_n\in B_b(D_m)$. By (\ref{fitz5}), Lemma \ref{lkjh} and (\ref{weis1}), we get
\begin{eqnarray}\label{cv2}
&&{\cal E}^0(\xi\chi_n,\phi)\nonumber\\
&=&\lim_{t\rightarrow
0}\frac{a^{\alpha}{\cal A}(d,-\alpha)}{t}\int_{D_m}\phi(x)\left[\int_0^t\int_{\overline{D}^c_m}\left(\xi(y)(1-\chi_n(y))
\int_{D_m}\frac{p^{D_m}(s,x,z)}{|z-y|^{d+\alpha}}dz\right)dyds\right]dx\nonumber\\
&=&\lim_{t\rightarrow
0}\frac{a^{\alpha}{\cal A}(d,-\alpha)}{t}\int_{D_m}\phi(x)\left[\int_0^t\int_{{D}^c_n}\left(\xi(y)(1-\chi_n(y))
\int_{D_m}\frac{p^{D_m}(s,x,z)}{|z-y|^{d+\alpha}}dz\right)dyds\right]dx\nonumber\\
&=&\lim_{t\rightarrow
0}\frac{a^{\alpha}{\cal A}(d,-\alpha)}{t}\int_0^t\int_{D_m}\phi(x)p^{D_m}_sF_n(x)dxds\nonumber\\
&=&a^{\alpha}{\cal A}(d,-\alpha)\int_{D_m}\phi(x)F_n(x)dx\nonumber\\
&=&a^{\alpha}{\cal A}(d,-\alpha)\int_{D_m}\int_{D^c_n}\frac{\xi(y)(1-\chi_n(y))}{|x-y|^{d+\alpha}}dy\phi(x) dx.
\end{eqnarray}
On the other hand, we have
\begin{eqnarray}\label{vb} & &{\cal
E}^0(\xi\chi_n,\phi)\nonumber\\&=&\int_{\mathbb{R}^d}\langle
\nabla(\xi\chi_n),\nabla\phi\rangle
dx-\int_{\mathbb{R}^d}\langle b,\nabla (\xi\chi_n)\rangle \phi dx\nonumber\\
&
&+\frac{a^{\alpha}{\cal A}(d,-\alpha)}{2}\int_{\mathbb{R}^d}\int_{\mathbb{R}^d}\frac{((\xi\chi_n)(x)-(\xi\chi_n)(y))(\phi(x)-\phi(y))}{|x-y|^{d+\alpha}}dxdy\nonumber\\
&=&{\cal E}^0(\xi,\phi)+a^{\alpha}{\cal A}(d,-\alpha)\int_{D_m}\int_{D^c_n}\frac{\xi(y)(1-\chi_n(y))}{|x-y|^{d+\alpha}}dy\phi(x) dx.
\end{eqnarray}
Thus, we obtain by (\ref{cv2}) and
(\ref{vb}) that ${\cal E}^0(\xi,\phi)=0$.

By (\ref{xw}), we have \begin{eqnarray*} &&\lim_{t\rightarrow
0}\int_Dw(x)\frac{w(x)-p^D_tw(x)}{t}dx\\
&\le&\lim_{t\rightarrow
0}\int_D|w|(x)\frac{E_x[\int_0^{t\wedge\tau}|f+cu|(X_s)ds]}{t}dx\\
&=&\lim_{t\rightarrow
0}\frac{1}{t}\int_0^t\left(p^D_s|f+cu|,|w|\right)ds\\
&=&(|f+cu|,|w|)\\
&<&\infty.
\end{eqnarray*}
Then, $w\in W^{1,2}_0(D)$ and hence $u=\xi+w\in W^{1,2}_{loc}(D)$.
For $\phi\in C^{\infty}_c(D)$, we have
\begin{eqnarray*} {\cal E}^0(w,\phi)&=&\lim_{t\rightarrow
0}\int_D\phi(x)\frac{w(x)-p^D_tw(x)}{t}dx\\
&=&\lim_{t\rightarrow
0}\int_D\phi(x)\frac{E_x[\int_0^{t\wedge\tau}(f+cu)(X_s)ds]}{t}dx\\
&=&\lim_{t\rightarrow
0}\frac{1}{t}\int_0^t\left(p^D_s(f+cu),\phi\right)ds\\
&=&(f+cu,\phi).
\end{eqnarray*}
Therefore,
$$
{\cal E}^0(u,\phi)={\cal E}^0(\xi+w,\phi)=(f+cu,\phi),
$$
which implies that (\ref{1}) holds.

\subsection{Uniqueness of solutions}

In this subsection, we will prove the uniqueness of solutions. To
this end, we will show that there exists $M>0$ such that if
$\|c^+\|_{L^{p\vee 1}}\le M$, then $v\equiv0$ is the unique
function in $B_b(\mathbb{R}^d)$ satisfying $v|_D\in
W^{1,2}_{loc}(D)\cap C(D)$ and
\begin{equation}\label{bl1}
\left\{\begin{array}{l}{\cal E}^0(v,\phi)=(cv,\phi),\ \ \forall \phi\in C^{\infty}_c(D),\\
v=0\ \ {\rm on}\ D^c.
\end{array}\right.\end{equation}

Suppose that $v\in B_b(\mathbb{R}^d)$ satisfying $v|_D\in
W^{1,2}_{loc}(D)\cap C(D)$ and (\ref{bl1}). Let
$\{D_n\}_{n\in\mathbb{N}}$ be a sequence of relatively compact
open subsets of $D$ such that $\overline{D}_n\subset D_{n+1}$ and
$D=\cup_{n=1}^{\infty}D_n$, and $\{\chi_n\}_{n\in\mathbb{N}}$  be
a sequence of functions in $C^{\infty}_c(D)$ such that $0\le
\chi_n\le 1$ and $\chi_n|_{D_n}=1$. We have $v\chi_n\in
W^{1,2}_0(D)$. Note that

\begin{equation}\label{san}
\int_{\mathbb{R}^d}\int_{\mathbb{R}^d}\frac{|(v(y)-v(x))(\chi_n(y)-\chi_n(x))|}{|x-y|^{d+\alpha}}dydx<\infty.
\end{equation}
Let $\beta>\beta_0$ (see (\ref{beta0})) and $\phi\in C^{\infty}_c(D)$. By (\ref{bl1}) and (\ref{san}), we get
\begin{eqnarray}\label{FOT}
& &{\cal
E}_{\beta}^0(v\chi_n,\phi)\nonumber\\&=&\int_{\mathbb{R}^d}\langle
\nabla(v\chi_n),\nabla\phi\rangle
dx-\int_{\mathbb{R}^d}\langle b,\nabla (v\chi_n)\rangle \phi dx\nonumber\\
&
&+\frac{a^{\alpha}{\cal A}(d,-\alpha)}{2}\int_{\mathbb{R}^d}\int_{\mathbb{R}^d}\frac{((v\chi_n)(x)-(v\chi_n)(y))(\phi(x)-\phi(y))}{|x-y|^{d+\alpha}}dxdy+(\beta,v\chi_n\phi)\nonumber\\
&=&{\cal E}^0(v,\chi_n\phi)-\int_{\mathbb{R}^d}(L\chi_n)v\phi
dx-2\int_{\mathbb{R}^d}\langle \nabla v,\nabla\chi_n\rangle\phi
dx\nonumber\\
&
&-a^{\alpha}{\cal A}(d,-\alpha)\int_{\mathbb{R}^d}\left[\int_{\mathbb{R}^d}\frac{(v(y)-v(x))(\chi_n(y)-\chi_n(x))}{|x-y|^{d+\alpha}}dy\right]\phi(x)dx+(\beta,v\chi_n\phi)\nonumber\\
&=&\left((c+\beta)v\chi_n-(L\chi_n)v-2\langle \nabla
v,\nabla\chi_n\rangle-a^{\alpha}{\cal A}(d,-\alpha)\int_{\mathbb{R}^d}\frac{(v(y)-v(\cdot))(\chi_n(y)-\chi_n(\cdot))}{|\cdot-y|^{d+\alpha}}dy,\phi\right)\nonumber\\
&:=&(\theta_n,\phi).
\end{eqnarray}

Let $n>m$ and $\phi\in C^{\infty}_c(D_m)$. By (\ref{FOT}), we get
\begin{eqnarray}\label{addsept}
& &(\theta_n,\phi)\nonumber\\
&=&{\cal E}^0(v,\phi)+a^{\alpha}{\cal A}(d,-\alpha)\int_{D_m}\int_{D^c_n}\frac{v(y)(1-\chi_n(y))}{|x-y|^{d+\alpha}}dy\phi(x) dx+(\beta,v\phi)\nonumber\\
&=&\left((c+\beta)v+a^{\alpha}{\cal A}(d,-\alpha)\int_{D\cap
D^c_n}\frac{v(y)(1-\chi_n(y))}{|\cdot-y|^{d+\alpha}}dy,\phi\right).
\end{eqnarray}
Since $\phi\in C^{\infty}_c(D_m)$ is arbitrary, by (\ref{addsept}), we find that for $n>m$,
$$
\theta_n(x)=(c(x)+\beta)v(x)+a^{\alpha}{\cal A}(d,-\alpha)\int_{D\cap
D^c_n}\frac{v(y)(1-\chi_n(y))}{|x-y|^{d+\alpha}}dy,\ \ x\in D_m.
$$
Hence
\begin{equation}\label{nm2}\theta_n\ {\rm converges\ to}\ (c+\beta)v\ {\rm uniformly\ on\ any\ compact\ subset\
of}\ D.
\end{equation}

Denote by
$((X_t)_{t\ge 0}, (P^{\beta}_x)_{x\in \mathbb{R}^d})$ the Markov
process associated with $({\cal E}_{\beta}^0,W^{1,2}(\mathbb{R}^d))$. For $n>m$, define
$$
A^{m,n}_t:=\int_0^{t\wedge\tau_{D_m}}\theta_n(X_s)ds\ \ {\rm and}\ \
c^{m,n}_t(x):=E^{\beta}_x[A^{m,n}_t],\ \ t\ge0,\ x\in D_m.
$$
By the joint continuity of $p(t,x,y)$ on $(0,\infty)\times
\mathbb{R}^d\times \mathbb{R}^d$, we know that the function $t\mapsto c^{m,n}_t(x)$ is
continuous for any $x\in D_m$. We have $c^{m,n}_t\in L^2(D_m;dx)$ for
$t\ge 0$ and
\begin{equation}\label{nm4}
c^{m,n}_{t+s}(x)=c^{m,n}_t(x)+p_t^{\beta,D_m}c^{m,n}_s(x),\ \ t,s\ge 0,
\end{equation}
where $(p_t^{\beta,D_m})_{t\ge 0}$ is the transition semigroup of the part
process $((X^{D_m}_t)_{t\ge 0},(P^{\beta}_x)_{x\in D_m})$.
By (\ref{FOT}), we get
\begin{eqnarray}\label{nm5}
\lim_{t\rightarrow 0}\frac{1}{t}E^{\beta}_{\phi\cdot dx}[A^{m,n}_t]&=&\lim_{t\rightarrow 0}\frac{1}{t}\int_0^t(p_s^{\beta,{D_m}}\theta_n,\phi)ds\nonumber\\
&=&(\theta_n,\phi)\nonumber\\
&=&{\cal
E}^0_{\beta}(v\chi_n,\phi),\ \ \ \ \forall \phi\in W^{1,2}_0(D_m).
\end{eqnarray}
Define
$$
\eta_{m,n}(x)=E^{\beta}_x[(v\chi_n)(X_{\tau_{D_m}})],\
\ x\in\mathbb{R}^d.
$$
We have
\begin{equation}\label{reduced}
\eta_{m,n}(x)=E^{\beta}_x[\eta_{m,n}(X_{t\wedge\tau_{D_m}})], \ \ t\ge 0,\ x\in D_m,
\end{equation} and $\eta_{m,n}(x)=v\chi_n(x)$ for
q.e.-$x\in D^c_m$. By \cite[Theorem 3.5.1]{Oshima}, we get
\begin{equation}\label{wanle}
{\cal E}^0_{\beta}(v\chi_n,\phi)={\cal
E}^0_{\beta}(v\chi_n-\eta_{m,n},\phi), \ \ \forall\phi\in W^{1,2}_0(D_m).
\end{equation}

Let $({T}^{\beta,D_m}_t)_{t\ge 0}$ be the $L^2$-semigroup associated
with $({\cal E}_{\beta}^0,W^{1,2}_0(D_m))$. Denote by $(\hat{T}^{\beta,D_m}_t)_{t\ge 0}$  the dual semigroup of $({T}^{\beta,D_m}_t)_{t\ge 0}$ on $L^2(D_m;dx)$. Define
\begin{equation}\label{eas}
\hat{S}^m_t:=\int_0^t\hat{T}^{\beta,D_m}_sds,\ \ t\ge 0.
\end{equation}
Similar to \cite[(1.5.5), page
39]{Fuku}, we can show that\begin{equation}\label{nm6} {\cal
E}_{\beta}^0(v,\hat{S}^m_t\rho)=(v,\rho-\hat{T}_t^{\beta,D_m}\rho),\ \ \forall
v\in W^{1,2}_0(D_m),\ \rho\in L^2(D_m;dx).
\end{equation}
Then, we obtain by (\ref{nm4}), (\ref{nm5}), (\ref{wanle}), (\ref{eas}) and (\ref{nm6}) that for $\phi\in
C^{\infty}_c(D_m)$ and $t,r>0$,
\begin{eqnarray*}
&&(c^{m,n}_t,\phi-\hat{T}^{\beta,D_m}_r\phi)\\&=&\lim_{s\rightarrow
0}\frac{1}{s}(c^{m,n}_t,\hat{S}^m_r\phi-\hat{T}^{\beta,D_m}_s\hat{S}^m_r\phi)\\
&=&\lim_{s\rightarrow
0}\frac{1}{s}(c^{m,n}_s,\hat{S}^m_r\phi-\hat{T}^{\beta,D_m}_t\hat{S}^m_r\phi)\\
&=&{\cal
E}^0_{\beta}(v\chi_n,\hat{S}^m_r\phi-\hat{T}^{\beta,D_m}_t\hat{S}^m_r\phi)\\
&=&{\cal
E}^0_{\beta}(v\chi_n-\eta_{m,n},\hat{S}^m_r\phi-\hat{T}^{\beta,D_m}_t\hat{S}^m_r\phi)\\
&=&(v\chi_n-\eta_{m,n},\phi-\hat{T}^{\beta,D_m}_t\phi-\hat{T}^{\beta,D_m}_r\phi+\hat{T}^{\beta,D_m}_{t+r}\phi)\\
&=&(v\chi_n-\eta_{m,n}-p^{\beta,D_m}_t(v\chi_n-\eta_{m,n}),\phi-\hat{T}^{\beta,D_m}_r\phi).
\end{eqnarray*}
Hence
$l^{m,n}_t:=(c^{m,n}_t-(v\chi_n-\eta_{m,n})+p^{\beta,D_m}_t(v\chi_n-\eta_{m,n}),\phi)$
satisfies the linear equation $l^{m,n}_t=l^{m,n}_{t+r}-l^{m,n}_r$.
By (\ref{nm5}) and (\ref{wanle}), we get $\lim_{t\rightarrow
0}l^{m,n}_t/t=0$. Then, $l^{m,n}_t=0$. Since $\phi\in
C^{\infty}_c(D_m)$ is arbitrary, we obtain by the continuity of
the function $t\mapsto p^{\beta,D_m}(t,x,y)$, which can be proved
similar to Lemma \ref{add1} (4), and the continuity of the
function $t\mapsto c^{m,n}_t(x)$ that for $dx$-a.e. $x\in D_m$,
$$
(v\chi_n-\eta_{m,n})(x)=E^{\beta}_x[(v\chi_n-\eta_{m,n})(X_{t\wedge\tau_{D_m}})]+E^{\beta}_x\left[\int_0^{t\wedge\tau_{D_m}}\theta_n(X_s)ds\right], \
\ \forall t\ge 0.
$$
By (\ref{reduced}), we obtain that for $dx$-a.e. $x\in D_m$,
\begin{equation}\label{zl}
(v\chi_n)(x)=E^{\beta}_x[(v\chi_n)(X_{t\wedge\tau_{D_m}})]+E^{\beta}_x\left[\int_0^{t\wedge\tau_{D_m}}\theta_n(X_s)ds\right], \
\ \forall t\ge 0.
\end{equation}

 Note that $v\in
B_b(\mathbb{R}^d)$ and $v=0$ on $D^c$. Letting
$n\rightarrow\infty$, we obtain by (\ref{nm2}) and (\ref{zl}) that
for $dx$-a.e. $x\in D_m$,
\begin{equation}\label{r2}
v(x)=E^{\beta}_x[v(X_{t\wedge\tau_{D_m}})]+E^{\beta}_x\left[\int_0^{t\wedge\tau_{D_m}}((c+\beta)v)(X_s)ds\right], \ \
\forall t\ge 0.
\end{equation}
Letting $m\rightarrow\infty$, we obtain by (\ref{r2}) that for
$dx$-a.e. $x\in D$,
\begin{equation}\label{chan}
v(x)=E^{\beta}_x[v(X_t)1_{\{\tau>t\}}]+E^{\beta}_x\left[\int_0^{t\wedge\tau}((c+\beta)v)(X_s)ds\right],
\ \ \forall t\ge 0.
\end{equation}

Define \begin{equation}\label{r1} {\cal
I}_t=v(X_t)1_{\{\tau>t\}}+\int_0^{t\wedge
\tau}((c+\beta)v)(X_s)ds.
\end{equation}
By (\ref{chan}), we find that $({\cal I}_t)_{t\ge 0}$ is a
martingale under $P^{\beta}_x$  for $dx$-a.e. $x\in D$. Define
$$
e_{\beta}(t):=e^{\int_0^t(c+\beta)(X_s)ds},\ \ t\ge0.
$$
The
integration by parts formula for semi-martingales implies that
$$
e_{\beta}(t){\cal I}_t-v(x)=\int_0^t{\cal
I}_sde_{\beta}(s)+\int_0^te_{\beta}(s)d{\cal I}_s.
$$By
(\ref{r1}), we get \begin{eqnarray*} && e_{\beta}(t){\cal
I}_t-\int_0^t{\cal
I}_sde_{\beta}(s)\nonumber\\&=&e_{\beta}(t)v(X_t)1_{\{\tau>t\}}+e_{\beta}(t)\int_0^{t\wedge
\tau}((c+\beta)v)(X_s)ds-\int_0^te_{\beta}(s)((c+\beta)v)(X_s)1_{\{\tau>
s\}}ds\nonumber\\
& &-\int_0^t((c+\beta)v)(X_w)1_{\{\tau\ge w\}}\left(\int_w^t
e_{\beta}(s)(c+\beta)(X_s)ds\right)dw\nonumber\\
&=&e_{\beta}(t)v(X_t)1_{\{\tau>t\}}\nonumber\\
&:=&{\cal J}_t.
\end{eqnarray*}
Hence $({\cal J}_t)_{t\ge 0}$ is a martingale under $P^{\beta}_x$
 for $dx$-a.e.
$x\in D$. Then, we have
\begin{eqnarray}\label{pp1}
v(x)&=&E^{\beta}_x[e_{\beta}(t)v(X_t)1_{\{\tau>t\}}]\nonumber\\
&=&E_x[e(t)v(X_t)1_{\{\tau>t\}}], \ \ dx-{\rm a.e.}\ x\in D.
\end{eqnarray}
By (\ref{14}), there exists $M>0$ such that if $\|c^+\|_{L^{p\vee
1}}\le M$ then
\begin{equation}\label{pp2}
\sup_{x\in D}E_x\left[e^{\int_0^{\tau}c^+(X_s)ds}\right]<\infty.
\end{equation}
Therefore, by letting $t\rightarrow\infty$, we obtain by
(\ref{pp1}), (\ref{pp2}) and the dominated convergence theorem
that $v(x)=0$ for $dx$-a.e. $x\in D$. Since $v|_D\in
 C(D)$, we obtain $v\equiv0$ on $\mathbb{R}^d$. The proof is complete.\hfill\fbox

\bigskip

{ \noindent {\bf\large Acknowledgments} \quad   This work was
supported by Natural Sciences and Engineering Research
Council of Canada. We thank the referee for the careful reading of our paper and all of the insightful
comments that greatly improved the presentation of the paper.}

\end{document}